\documentclass[10pt]{article}
\usepackage{amsmath}
\usepackage{amssymb}
\newcommand{\ppip}{\partial\Pi_+}
\newcommand{\pip}{\Pi_+}
\newcommand{\hht}[1]{{\widehat #1}}
\newcommand{\phiz}{\phi_z}
\newcommand{\eisx}{e^{-isx}}
\newcommand{\episx}{e^{isx}}
\newcommand{\eis}[1]{e^{-is{#1}}}
\newcommand{\R}{{\mathbb R}}
\newcommand{\rbar}{{\overline {\mathbb R}}}
\newcommand{\fn}{\!:\!}
\newcommand{\real}{\!:\!\R\to\R}

\newcommand{\lsum}{\sum\limits}
\newcommand{\llim}{\lim\limits}
\newcommand{\lint}{\int\limits}
\newcommand{\linf}{\inf\limits}
\newcommand{\lsup}{\sup\limits}
\newcommand{\qed}{\mbox{$\quad\blacksquare$}}
\newcommand{\hk}{{\cal HK}}
\newcommand{\bv}{{\cal BV}}
\newcommand{\meas}{{\cal M}}

\newcommand{\fhat}{{\widehat f}}

\newtheorem{theorem}{Theorem}
\newtheorem{lemma}[theorem]{Lemma}
\newtheorem{prop}[theorem]{Proposition}
\newtheorem{corollary}[theorem]{Corollary}
\newtheorem{remark}[theorem]{Remark}
\newtheorem{example}[theorem]{Example}
\newtheorem{defn}[theorem]{Definition}
\begin{document}
\begin{center}
{\large\bf Henstock--Kurzweil Fourier transforms
}
\vskip.25in
Erik Talvila\footnote{Research partially supported by NSERC.}\\ [2mm]
{\footnotesize
Department of Mathematical and Statistical Sciences\\
University of Alberta\\
Edmonton AB Canada T6G 2E2\\
etalvila@math.ualberta.ca}
\end{center}

{\footnotesize
\noindent
{\bf Abstract:} The Fourier transform is considered as a Henstock--Kurzweil integral.
Sufficient conditions are given for the existence of the Fourier transform and
necessary and sufficient conditions are given for it to be continuous. 
The Riemann--Lebesgue lemma fails:  Henstock--Kurzweil Fourier transforms
can have  arbitrarily large point-wise growth.  Convolution and inversion theorems
are established.  An appendix gives sufficient conditions for interchanging
repeated Henstock--Kurzweil integrals and gives an estimate on the integral of a product.
\\
{\bf 2000 subject classification:} 42A38, 26A39
}
\vskip.25in

\section{Introduction}
If $f\real$ then its Fourier transform at $s\in\R$ is defined as
$\fhat(s)=\int_{-\infty}^\infty \eisx f(x)\,dx$.  The inverse transform
is $f\!\!\text{{\rm{\LARGE \v{}}}}\!\!(s)=(2\pi)^{-1}\int_{-\infty}^\infty e^{isx} f(x)\,dx$.
In this paper we consider Fourier transforms as Henstock--Kurzweil integrals.
This is an integral equivalent to the Denjoy and Perron integrals but 
with a definition in terms of Riemann sums.
We let $\hk_A$ be
the Henstock--Kurzweil integrable functions over a set $A\subset\R$, 
dropping the subscript
when $A=\R$.  (The symbol $\subset$ allows set equality.)
Then $\hk$ properly contains the union of $L^1$ and the Cauchy-Lebesgue
integrable functions (i.e., improper Lebesgue integrals).  The main points of
$\hk$ integration that we use can be found in \cite{bartle} and \cite{swartz}.
Several of our results depend on being able to reverse the order of repeated
integrals.  In the Lebesgue theory this is usually justified with Fubini's Theorem.
For $\hk$ integrals, necessary and sufficient conditions were given in 
\cite{talviladiff}.  Lemma~\ref{lemma3} in the Appendix gives sufficient conditions
that are readily applicable to the cases at hand.  Also in the Appendix are
some conditions for convergence  of rapidly oscillatory integrals 
(Lemma~\ref{lemma1})
and 
an estimate of the integral of a product (Lemma~\ref{lemma2}).

There is a substantial body of theory relating to Fourier transforms when
they are considered as Lebesgue integrals.  Necessary and sufficient for
existence of $\fhat$ on $\R$ is that $f\in L^1$.  This is because the multipliers
for $L^1$ are the (essentially) bounded measurable functions and $|e^{\pm isx}|\leq 1$.  The
multipliers for $\hk$ are the functions of (essentially) bounded variation.  As
$x\mapsto e^{-isx}$ is not of bounded variation, 
except for $s=0$, we do not have an elegant
existence theorem for $\hk$ Fourier integrals.  Various existence conditions are given
in Proposition~\ref{prop1}.  Example~\ref{example1}(f)
gives a function whose Fourier transform
diverges on a countable set.   For $L^1$ convergence, $\fhat$ is uniformly continuous
with limit $0$ at infinity (the
Riemann--Lebesgue
lemma). 
We show below (Example~\ref{example1}(e)) that the Riemann--Lebesgue lemma fails dramatically in $\hk$:
$\fhat$ can have arbitrarily large point-wise growth.
And, $\fhat$ need not be continuous.   Continuity of $\fhat$ is equivalent to 
quasi-uniform convergence
(Theorem~\ref{quc}).
Some sufficient conditions for continuity of $\fhat$ appear in 
Proposition~\ref{propcts}.
  Although $\fhat$ need not be continuous, when it exists at the
endpoints of a compact interval, it exists almost everywhere on that interval and
is integrable over that interval.  See Proposition~\ref{b1}.
As in the $L^1$ theory, we have
linearity, symmetry, conjugation, translation, modulation,
dilation, etc.
See formulas (2)--(9) in \cite[p.~117]{erdelyi} and \cite[p.~9]{benedetto}.
We draw attention to the differentiation of Fourier transforms (Proposition~\ref{freq}) and
transforms of derivatives (Proposition~\ref{time}).  One of the most important properties
of Fourier transforms is their interaction with convolutions.
Propositions~\ref{c1}, \ref{convexist}, \ref{c3}, \ref{c4} and \ref{c5} contain various results
on existence of convolutions; estimates using the variation, $L^1$ norm and Alexiewicz norm;
and the transform and inverse transform of convolutions.
Proposition~\ref{propinv} gives a Parseval relation.
An inversion theorem is obtained using a summability kernel (Theorem~\ref{inversion}).  A
uniqueness theorem follows as a corollary.  The paper concludes with
an example of a function $f$ for which $\fhat$ exists on $\R$ but 
$\fhat\text{{\rm{\LARGE \v{}}}}$ exists nowhere.

As Henstock--Kurzweil
integrals allow conditional convergence, they make an ideal setting for the
Fourier transform.  We remark that many of the Fourier integrals appearing in
tables such as \cite{erdelyi} diverge as Lebesgue integrals but converge as
improper Riemann integrals.  Thus, they exist as $\hk$
integrals.

We use the following notation.
Let $A\subset\R$ and $f$ a real-valued function
on $A$.  
The functions of bounded variation over $A$ are denoted
$\bv_A$ and the variation of function $f$ over set $A$ is
$V_Af$.  We say a set is in $\bv$ if its characteristic function is
in $\bv$.  
All our results are stated for real-valued functions but the extension
to complex-valued functions is immediate.  Note that for complex-valued functions, 
the variation of the real part and the
variation of the imaginary part are added. The Alexiewicz norm of $f\in\hk_A$ is
$\|f\|_A=\sup_{I\subset A}|\int_I f|$,  the supremum being taken over all
intervals $I\subset A$.  
For each of these definitions, the label $A$  is  
omitted when $A=\R$ or it is obvious which set is $A$.  
Whereas indefinite Lebesgue integrals are absolutely continuous ($AC$), indefinite
Henstock--Kurzweil integrals are
$ACG_*$.  See \cite{saks} for the definition of $ACG_*$ and the related
space $AC_*$.
Finally, a convergence theorem that we use throughout is:
\begin{theorem}
Let $f$ and $g_n$ $(n\in\mathbb{N})$ be real-valued 
functions on $[a,b]$.  If $f\in\hk$,
$Vg_n\leq M$ for all $n\in\mathbb{N}$, and $g_n\to g$ as $n\to\infty$
then $\int_a^b fg_n\to\int_a^b fg$ as $n\to\infty$.
\end{theorem}
The theorem holds for $[a,b]\subset\rbar$, where $\rbar=\R\cup\{\pm\infty\}$ is
the extended
real line.  
For a proof see \cite{talvilarae}.

\section{Basic properties}
We first tackle the problem of existence. If $f\real$ then $\fhat$ exists as
a Lebesgue integral on $\R$ if and only if $f\in L^1$.  This
follows from the fact that $|e^{\pm isx}|\leq 1$ for all $s,x\in\R$ and
the multipliers for $L^1$ are the bounded measurable functions.
No such simple necessary
and sufficient conditions are known for existence of $\hk$ Fourier integrals.
However, we do have the following results.

\begin{prop}\label{prop1}
Let $f\real$. 
\begin{enumerate}
\item[{\rm (a)}] 
In order for $\fhat$ to exist at some $s\in\R$ it is necessary that $f\in\hk_{loc}$.
\item[{\rm (b)}]
If $f\in\hk_{loc}$ then $\fhat$ exists on $\R$ if $|f|$ is integrable in a 
neighbourhood of infinity or if $f$ is 
of bounded variation in a neighbourhood of infinity with limit $0$ at infinity.
\item[{\rm (c)}]
Let $f\in\hk$.  Define $F_1(x)=\int_x^{\infty}f$ and 
$F_2(x)=\int_{-\infty}^xf$.  Then $\fhat$ exists at 
$s\in\R$ if and only if  both the integrals $\int_0^\infty\eisx F_1(x)\,dx$ and 
$\int_{-\infty}^0\eisx F_2(x)\,dx$ exist.
\end{enumerate}
\end{prop}

\noindent
{\it Proof:} (a) For each $s\in\R$, the function $x\mapsto \episx$ is 
of bounded variation on 
any compact interval.\\

(b) This follows from the Chartier--Dirichlet convergence test.  See 
\cite{bartle}.\\

(c) Let $T>0$.  Integrate by parts to obtain
$$
\int_{0}^T\eisx f(x)\,dx=F_1(0)-F_1(T)\eis{T} -is\int_0^T\eisx F_1(x)\,dx.
$$
Since $F_1$ is continuous with limit $0$ at infinity, $\int_{0}^{\infty}\eisx f(x)\,dx$
exists if and only if $\int_0^{\infty}\eisx F_1(x)\,dx$ exists.  The other part of
the proof is similar. \qed\\

Although $F_1$ is continuous with limit $0$ at infinity, it need not be of
bounded variation.  So, $f\in\hk$ does not imply the existence of $\fhat$.  See
Example~\ref{example1}(c)  below.
Notice that part (b) (with $\hk_{loc}$ replaced by $L^1_{loc}$)  and part
(c) are false for $L^1$ convergence of $\fhat$.

Titchmarsh \cite{titchmarsh} gives several sufficient conditions for existence
of conditionally convergent Fourier integrals (\S1.10--1.12).  However, these
all require that $f\in L^1_{loc}$.

When $f\in L^1$ and $s,h\in\R$, $\fhat(s+h)=\int_{-\infty}^{\infty}e^{-i(s+h)x}f(x)\,dx$.
By dominated convergence this tends to $\fhat(s)$ as $s\to h$.  So, $\fhat$ is
uniformly continuous on $\R$.  When $\fhat$ exists in $\hk$ in a neighbourhood
of $s$, the function $x\mapsto \eisx f(x)$ is in $\hk$ but the factor 
$e^{-ihx}$ is not
of bounded variation on $\R$ except for $h=0$.  In general we cannot take the
limit $h\to 0$ under the integral sign and $\fhat$ need not be continuous.
And, for $f\in L^1$ and $s\not=0$, the 
change of variables $x\mapsto x+\pi/s$ gives
$\fhat(s)=(1/2)\int_{-\infty}^\infty \eisx\left[f(x)-f(x+\pi/s)\right]\,dx$.
Writing $f_y(x)=f(x+y)$ for $x,y\in\R$, we have $|\fhat(s)|\leq (1/2)\|f-f_{\pi/s}\|_1$.
Continuity of $f$ in the $L^1$ norm now yields the Riemann--Lebesgue lemma:
$\fhat(s)\to 0$ as $|s|\to\infty$.  It is true that if $f\in\hk$ then $f$ is continuous
in the Alexiewicz norm \cite{talvilacalex}.  However, since the variation of $x\mapsto \eisx$ is not 
uniformly bounded  as $|s|\to\infty$, existence of $\fhat$ does not let us conclude
that $\fhat$ tends to $0$ at infinity.

The following examples show some of the differences between $L^1$ and $\hk$ 
Fourier transforms.

\begin{example}\label{example1}
{\rm
The transforms (a)--(d) appear in \cite{erdelyi}.  Convergence in (a) is 
by Lemma~\ref{lemma1}, (b) is similar, after integrating by parts,
and (c) and (f) are Frullani integrals.\\

\noindent
(a)\quad If $f(x)={\rm sgn}(x)|x|^{-1/2}$ then $f$ is not in $\hk$ or in any $L^p$ space 
($1\leq p\leq\infty$)
and yet $\fhat(s)=\sqrt{2\pi}\,{\rm sgn}(s)|s|^{-1/2}$ for $s\not=0$.  Notice that,
even though $f$ is odd, $\fhat$ does not
exist at $0$ since $\hk$ convergence does not allow principal value integrals.  \\

\noindent
(b)\quad 
Let $g(x)=e^{ix^2}$.  Then $\hht{g}(s)=\sqrt{\pi}\,e^{i(\pi-s^2)/4}$.  In this 
example,
$\hht{g}$ is not of bounded variation at infinity, nor does $\hht{g}$ tend to $0$ at
infinity, nor is $\hht{g}$ uniformly continuous on $\R$.  The same can of course
be said for $g$.\\

\noindent
(c) \quad Let $h(x)=\sin(ax)/|x|$.  Then $\hht{h}(s)=i\log|(s-a)/(s+a)|$ for $s\not=a$.\\

\noindent
(d)\quad Let $k(x)=x/(x^2+1)$.  Then $\hht{k}(s)=-i\pi\,{\rm sgn}(s)e^{-|s|}$ for
$s\not=0$.  Note that $\hht{k}$ does not exist at $0$, even though its principal
value is $0$.\\

\noindent
(e)\quad Fourier transforms
in $\hk$ can have arbitrarily large point-wise growth.
Given any sequence $\{a_n\}$ of positive real numbers, 
there is a continuous function
$f$ on $\R$ such that $\fhat$ exists 
on $\R$ and $\fhat(n)\geq a_n$ for all $n\geq 1$ \cite{talvilarapid}.
\\

\noindent
(f) \quad Let $\{a_n\}$ and $\{b_n\}$ be  sequences in $\R$.
Define $f(x)=\sum_{n=1}^\infty a_n\sin(b_nx)/|x|$ for
$x\not=0$ and $f(0)=0$.  Assume that $a_n>0$,  $\sum a_n<\infty$ and
$\sum a_n|b_n|<\infty$.  Then $f$ is continuous on $\R$, except at the origin,
where it has a finite jump discontinuity.  Suppose $s$ is not in
the closure of $\{-b_n,b_n\}_{n\in\mathbb{N}}$.  Then
\begin{eqnarray}
\fhat(s) & = & \lsum_{n=1}^\infty a_n \int_{-\infty}^\infty e^{-isx}\sin(b_nx)\,
\frac{dx}{|x|}\label{f1}\\
 & = & i\lsum_{n=1}^\infty a_n \int_{0}^\infty \left(\cos\left[(s+b_n)x\right]
 -\cos\left[(s-b_n)x\right]\right)\frac{dx}{x}\label{f2}\\
 & = & i\lsum_{n=1}^\infty a_n\log\left|\frac{s-b_n}{s+b_n}\right|.\label{f3}
\end{eqnarray}
The reversal of summation and integration in \eqref{f1} is justified using
Corollary~7 in \cite{talviladiff}.  Hence, $\fhat$ exists on $\R$, except
perhaps on the closure of $\{-b_n,b_n\}_{n\in\mathbb{N}}$.  Note that
$\fhat(0)=0$. 

We will now show $\fhat$ diverges at
each $b_k$ with $a_kb_k\not=0$.  Let $T_1, T_2>0$ and consider
\begin{eqnarray}
\lefteqn{
\int_{-T_1}^{T_2} e^{-ib_kx}\lsum_{n=1}^\infty a_n\sin(b_nx)
\frac{dx}{|x|}}\notag\\
& = & \lsum_{n=1}^\infty a_n\!\!\int_{-T_1}^{T_2} e^{-ib_kx}\sin(b_nx)\frac{dx}{|x|}
\label{f4}\\
 & = & \!\!\lsum_{n=1}^\infty\! a_n\!\!\int_{\!-T_1}^{T_2}\!\!\frac{
	 \sin\left[(b_k+b_n)x\right]
	  \!-\!\sin\left[(b_k-b_n)x\right]}{2}-i\sin(b_kx)\sin(b_nx)\frac{dx}{|x|}.
	  \label{f5}
\end{eqnarray}
In \eqref{f4}, convergence of $\sum a_n|b_n|$ permits reversal of summation
and integration.  The real part of \eqref{f5}
converges for all $k\geq 1$, uniformly for $T_1, T_2\geq 0$.  Hence, the
real part of $\fhat$ exists on $\R$.  The $k^{th}$ summand of the imaginary
part of \eqref{f5} is
$$
-a_k\int_{-T_1}^{T_2}\sin^2(b_kx)\frac{dx}{|x|}=
-a_k\int_{-T_1|b_k|}^{T_2|b_k|}\sin^2x\frac{dx}{|x|}.
$$
This diverges as $T_1,T_2\to\infty$.  Hence, $\fhat(b_k)$ does not exist.

If $\{-b_n,b_n\}_{n\in\mathbb{N}}$ has no limit points then we have an example of a function
whose Fourier transform exists everywhere except on a countable set.

Now suppose $s\not\in\{-b_n,b_n\}_{n\in\mathbb{N}}$ but $s$
is a limit point of $\{-b_n,b_n\}_{n\in\mathbb{N}}$.  
As noticed above,
the real part of $\fhat(s)$ exists.  And,
$|\int_{-1}^1\sin(sx)\sin(b_nx)\,dx/|x||\leq 2|s|$.  So, $\fhat(s)$ exists
if and only if
$$
\llim_{T\to\infty}\lsum_{n=1}^\infty a_n\!\!\int_{1}^{T}
\left(\cos\left[(s+b_n)x\right]
 -\cos\left[(s-b_n)x\right]\right)\frac{dx}{x}
$$
exists.  Suppose $s\not=0$ and $T>1$.
If $|s-b_n|T>1$ and $|s-b_n|<1$ then
\begin{eqnarray}
\left|\int_1^T\cos(|s-b_n|x)\,\frac{dx}{x}\right| & = &
\left|\int_{|s-b_n|}^{|s-b_n|T}\cos x\,\frac{dx}{x}\right|\notag\\
 & = & \left|\int_{|s-b_n|}^1\cos x\,\frac{dx}{x} + 
 \int_{1}^{|s-b_n|T}\cos x\,\frac{dx}{x}\right|\label{f5.5}\\
 & \leq & \log\left(1/|s-b_n|\right)+c.\notag
\end{eqnarray}
The constant $c$ is equal to the supremum of $|\int_1^t \cos x\,dx/x|$ over
$t>1$.  When $|s-b_n|T\leq 1$, we have
\begin{eqnarray*}
\left|\int_1^T\cos(|s-b_n|x)\,\frac{dx}{x}\right| & = &
\left|\int_{|s-b_n|}^{|s-b_n|T}\cos x\,\frac{dx}{x}\right|\\
& \leq & \log T\\
& \leq & \log\left(1/|s-b_n|\right).
\end{eqnarray*}
The case for $|s+b_n|T$ is similar.  It follows that the condition
\begin{equation}
\lsum_{n=1}^\infty a_n\left|\log\left|\frac{s-b_n}{s+b_n}\right|\right|<\infty
\label{f6}
\end{equation}
is sufficient for existence of $\fhat(s)$.

If $1/T<|s-b_n|< 1$ then, as in \eqref{f5.5},
$$
\int_1^T\cos(|s-b_n|x)\,\frac{dx}{x}\geq 
\cos(1)\log\left(1/|s-b_n|\right)-c.
$$
Therefore, 
\begin{eqnarray*}
\lefteqn{\!\!\!\!\!\!\!\!\!\!\!\!\!\!\!\!
\lsum_{1/T<|s-b_n|< 1} a_n\int_1^T\cos(|s-b_n|x)\,\frac{dx}{x}}\\ 
 & \geq & \lsum_{1/T<|s-b_n|< 1}a_n\left[\cos(1)\log\left(1/|s-b_n|\right)-c\right].
 \end{eqnarray*}
Let $T\to\infty$, then condition \eqref{f6} is also
necessary for existence of $\fhat(s)$.  Hence, it is possible for $\fhat$
to exist at a finite number of limit points of 
$\{-b_n,b_n\}_{n\in\mathbb{N}}.$

Finally, enumerate the rational numbers in $[0,1]$ by $b_1=0, b_2=1/1, b_3=1/2,
b_4=1/3, b_5=2/3, b_6=3/3, b_7=1/4$, etc.
Let $A_m>0$ such that
$\sum mA_m<\infty$.  Put $a_1=0$ and define $a_n=A_m$ for the $m$
consecutive values of $n$ such that $b_n=l/m$ for some
$1\leq l\leq m$.  Let $s\in[-1,1]\setminus{\mathbb Q}$ and let
$\bar{s}$ be the distance to the nearest rational number.  Then
$$
\lsum_{n=1}^\infty a_n\left|\log\left|\frac{s-b_n}{s+b_n}\right|\right|
=
\lsum_{m=1}^\infty A_m\lsum_{l=1}^m\left|\log\left|\frac{s-l/m}{s+l/m}\right|
\right|
\leq
\log(2/\bar{s})\lsum_{m=1}^\infty mA_m.
$$
This furnishes an example of a function whose Fourier transform exists
on $\R$ except for the rational numbers in $[-1,1]$.\qed
}
\end{example}

Examples~\ref{example1}(a), (c), (d) and (f)  show that
$\fhat$  need not be
continuous.   However, continuity of $\fhat$ is equivalent to quasi-uniform
continuity.

\begin{defn}[Quasi-uniform continuity]
Let $f\fn\R^2\to\R$.  If $F(x):=\int_{-\infty}^\infty f(x,y)\,dy$ exists in a 
neighbourhood
of $x_0\in\R$ then $F$ is quasi-uniformly continuous at $x_0$  if
for all $\epsilon>0$ and $M>0$ there exist $m=m(x_0,\epsilon,M)\geq
M$ and
$\delta=\delta(x_0,\epsilon,M)>0$ such that if $|x-x_0|<\delta$ then 
$|\int_{|y|>m}
f(x,y)\,dy|<\epsilon$.
\end{defn}
This is a modification of a similar definition for series, originally introduced by
Dini.  See \cite{bromwich}, page~140.

\begin{theorem}\label{quc}
Let $f\real$.
Then,
$\fhat$ is continuous at $s_0\in\R$ if and only if $\fhat$ is quasi-uniformly continuous 
at $s_0$.
\end{theorem}

\noindent
{\bf Proof:} For $m>0$, let $F_m(s)=\int_{-m}^m \eisx f(x)\,dx$.  Let $h\in\R$. Then,
$F_m(s+h)-F_m(s)=\int_{-m}^m[e^{-ihx}-1] \eisx f(x)\,dx$.  Note that either assumption
implies $x\mapsto \eisx f(x)$ is in $\hk_{loc}$ for each  $s\in\R$.
And, $V_{[-m,m]}[x\mapsto e^{-ihx}-1]\leq 
4m|h|$.  Taking the limit $h\to 0$ inside the above integral now shows $F_m$ is continuous
on $\R$ for each $m>0$.

Suppose $\fhat$ is quasi-uniformly continuous at $s_0\in\R$.  Given $\epsilon>0$, take
$M>0$ such that $|\int_{|x|>t}e^{is_0x}f(x)\,dx|<\epsilon$ for all $t>M$.  From 
quasi-uniform continuity, we have $m>M$ and $\delta>0$.  Then, for $|s-s_0|<\delta$,
\begin{eqnarray*}
\left|\,\lint_{|x|>m}\!\!\!\!\left[\eisx - e^{-is_0x}\right]f(x)\,dx\right| & 
\leq & 
\left|\,\lint_{|x|>m}\!\!\!\!\eisx f(x)\,dx\right| + \left|\,\lint_{|x|>m}\!\!\!\!
 e^{-is_0x}f(x)\,dx\right|\\
 & \leq & 2\epsilon.
\end{eqnarray*}
It follows that $\fhat$ is continuous at $s_0$.

Suppose $\fhat$ is continuous at $s_0$ and we are given $\epsilon >0$ and $M>0$.
Since $\fhat$ exists at $s_0$, there is $N=N(s_0,\epsilon)>0$ such that
$|\int_{|x|>m} e^{-is_0x}f(x)\,dx|<\epsilon$ whenever $m>N$.  Continuity
of $\fhat$ at $s_0$ implies the existence of $\xi=\xi(s_0,\epsilon)>0$ such that
$|\fhat(s)-\fhat(s_0)|<\epsilon$ when $|s-s_0|<\xi$.
And, $F_m$ is continuous on $\R$.  Hence, there exists 
$\eta=\eta(s_0,\epsilon,m)>0$
such that when $|s-s_0|<\eta$ we have
$|F_m(s)-F_m(s_0)|<\epsilon$.  Let $m=\max(M,N)$ and
$\delta=\min(\xi,\eta)$.  Then for
$|s-s_0|<\delta$ we have
\begin{eqnarray*}
\left|\,\lint_{|x|>m}\!\!\!\!\eisx f(x)\,dx\right| & \!\leq\! & \!\left|\,
\lint_{|x|>m}\!\!\!\!e^{is_0x} 
f(x)\,dx\right|
\!+\! \left|\fhat(s)-\fhat(s_0)\right|
\!+\! \left|F_m(s)\!-\!F_m(s_0)\right|\\
 & < & 3\epsilon.
\end{eqnarray*}
And, $\fhat$ is quasi-uniformly continuous at $s_0$.  \qed

We now present two sufficient conditions for a Fourier transform to be continuous.
The first is in the spirit of the Chartier--Dirichlet convergence test 
and the second is in the
spirit of the Abel convergence test.
For simplicity, the results are stated for functions on $[0,\infty)$.  The general case
follows easily.
\begin{prop}\label{propcts}
Let $g$ and  $h$ be real-valued functions on $[0,\infty)$ where $g\in\bv$
and $h\in\hk_{{\rm loc}}$.  Define $f=gh$.
\begin{enumerate}
\item[{\rm (a)}] Suppose 
there are positive constants $M$, $\delta$ and $K$ such that, if $|s-s_0|<\delta$
and $M_1, M_2>M$
then $|\int_{M_1}^{M_2}\eisx h(x)\,dx|< K$.  If $g(x)\to 0$ 
as $x\to\infty$ then $\fhat$ is continuous at $s_0$.
\item[{\rm (b)}]  Let $H_s(x)=\int_0^xe^{-ist}f(t)\,dt$. 
If $\hht{h}$ is continuous at $s_0$ and there are $\delta, K>0$ such that
for all $x>0$ and $|s-s_0|<\delta$ we have $|H_s(x)|\leq K$
then $\fhat$ is continuous at $s_0$. 
\end{enumerate}
\end{prop}

\noindent
{\bf Proof:}
Write $\phi_s(x)=\eisx h(x)$.
With no loss of generality, $g(\infty)=0$.

For (a), let $|s-s_0|<\delta$ and $M_1, M_2 > M$.
Using Lemma~\ref{lemma2},
\begin{eqnarray}
\left|\int_{M_1}^{M_2}\eisx f(x)\,dx\right| & \leq & 
\left|\int_{M_1}^{M_2}\phi_s(x)\,dx\right|\linf_{[M_1,M_2]}|g|+\|\phi_s\|_{[M_1,M_2]}
V_{[M_1,M_2]}g\notag\\
 & \leq & K\left[\linf_{[M_1,M_2]}|g|+V_{[M,\infty]}g\right]\label{f7}\\
 & \to & 0 \quad\text{ as } M\to \infty.\notag
 \end{eqnarray}
Therefore, $\fhat$ exists in a neighbourhood of $s_0$. Taking the limit
$M_2\to\infty$ in \eqref{f7} shows that
$\fhat$ is quasi-uniformly continuous and hence continuous.

For (b), since $g\in\bv$ we have $\lim_{x\to\infty}g(x)=c\in\R$.  Writing
$f=h(g-c)+ch$ we need only consider 
$|\int_M^\infty \phi_s(g-c)|\leq \|\phi_s\|_{[M,\infty)}V_{[M,\infty)}(g-c)$.
By our assumption, $\|\phi_s\|_{[M,\infty)}\leq 2K$ for $|s-s_0|<\delta$.
And, $V_{[M,\infty)}(g-c)\to 0$ as $M\to\infty$.\qed

Although $\fhat$ need not be continuous, when it exists  at the endpoints of a
compact interval it is integrable over the interval.

\begin{prop}\label{b1}
Let $[a,b]$ be a compact interval.  If $\fhat$ exists at $a$ and $b$ then 
$\fhat$ exists almost everywhere on $(a,b)$, $\fhat$ is integrable over
$(a,b)$ and $\int_a^b\fhat=i\int_{-\infty}^\infty f(x)[e^{-ibx}-e^{-iax}]\,dx/x$.
\end{prop}

\noindent
{\bf Proof:}
The integral $I:=i\int_{-\infty}^\infty f(x)\left[e^{-ibx}-e^{-iax}\right]\,
\frac{dx}{x}$ exists since $x\mapsto f(x)e^{-ibx}/x$ and $x\mapsto f(x)e^{-iax}/x$
are integrable over $\R\setminus(-1,1)$ and 
$x\mapsto\left[e^{-ibx}-e^{-iax}\right]/x$ is of bounded variation on
$[-1,1]$.  And,
$$
I  =  \int_{-\infty}^\infty f(x)e^{-ibx}\int_a^b e^{-i(s-b)x}ds\,dx
  =  \int_a^b\int_{-\infty}^\infty f(x)e^{-isx}dx\,ds
  =  \int_a^b \fhat.
$$
Hence, $\fhat$ exists almost everywhere on $(a,b)$ and is integrable over
$(a,b)$.
Lemma~\ref{lemma3}(a) justifies the reversal of $x$ and $s$ integration.\qed

The usual algebraic 
properties of linearity, symmetry, conjugation, translation, modulation,
dilation, etc., familiar from the $L^1$ theory, continue to hold for $\hk$ Fourier
transforms.  See formulas (2)--(9) in \cite[p.~117]{erdelyi} and \cite[p.~9]{benedetto}.
The proofs are elementary.  There are also differentiation results analogous to the 
$L^1$ case (pages  117 and 17, respectively, of the previous references).

\begin{prop}[Frequency differentiation]\label{freq}
Suppose $\fhat$ exists on  the compact interval $[\alpha,\beta]$.  
Define $g(x)=xf(x)$ and suppose $g\in \hk$.
Then
${\hat f}\,'=-i\hht{g}$
almost everywhere on $(\alpha,\beta)$.
In particular,
${\hat f}\,'(s)=-i\hht{g}(s)$ for all $s\in(\alpha,\beta)$ such that
$\frac{d}{ds}\int_{\alpha}^s\hht{g}=\hht{g}(s)$.
\end{prop}

\noindent
{\bf Proof:} The necessary and sufficient condition that allows
differentiation under the integral,
${\hat f}\,'(s)=-i\int_{-\infty}^\infty \eisx xf(x)\,dx$, for 
almost all $s\in(\alpha,\beta)$ is
that
\begin{equation}
\int_{-\infty}^\infty \int_a^b \eisx xf(x)\,ds\,dx=\int_a^b\int_{-\infty}^{\infty}\eisx xf(x)\,dx\,ds\label{b2}
\end{equation}
for all $[a,b]\subset[\alpha,\beta]$.  See \cite[Theorem 4]{talviladiff}.  
We have $g\in\hk$, $|\eisx|
\leq 1$ and $V_I[x\mapsto \eisx]\leq 2|I||s|$ for a compact interval  
$I\subset\R$.  
The left member of \eqref{b2} is $i[\fhat(b)-\fhat(a)]$.
Hence, by Lemma~\ref{lemma3}(a), 
\eqref{b2} holds, and ${\hat f}\,'(s)=-i\hht{g}(s)$ 
for almost all $s\in(\alpha,\beta)$.
Examining the proof of \cite[Theorem 4]{talviladiff}, we see that we get equality ${\hat f}\,'(s)=-i\hht{g}(s)$
when $\frac{d}{ds}\int_{\alpha}^s\hht{g}=\hht{g}(s)$.
\qed 

There are similar results for $n$-fold differentiation
when the function $x\mapsto x^nf(x)$ is in $\hk$ for a positive integer $n$.

\begin{prop}[Time differentiation]\label{time}
\begin{enumerate}
\item[{\rm (a)}] If $f\in ACG_*(\R)$ and $f(x)\to 0$ as $|x|\to \infty$ 
then for each $s\not=0$,
both $\fhat(s)$ and  $\hht{f'}(s)$ fail to exist or $\hht{f'}(s)=is\fhat(s)$.
\item[{\rm (b)}] Suppose $f\in ACG_*(\R)$ and $f, f'\in\hk$.  Then for each $s\not=0$, either
both $\fhat(s)$ and  $\hht{f'}(s)$ fail to exist or $\hht{f'}(s)=is\fhat(s)$.
\end{enumerate}
\end{prop}

\noindent
{\bf Proof:} 
(a) Let $M_1,M_2>0$.  Integrate by parts to get
$$
\int_{-M_1}^{M_2}\eisx f'(x)\,dx = e^{-isM_2}f(M_2)-e^{isM_1}f(-M_1)
+is\int_{-M_1}^{M_2}\eisx f(x)\,dx.
$$
Now take the limits $M_1, M_2\to\infty$.

(b) Consider $\int_x^{\infty}f'=\int_M^{\infty}f'+f(M)-f(x)$ for $x,M\in\R$.
Since $f'\in\hk$, the limits as $|x|\to\infty$ exist.  Hence, $f$ has a limit
at infinity.  But, $f\in\hk$ so this limit must be $0$ and we have reduction to case~(a).
\qed

\section{Convolution}
If $f$ and $g$ are real-valued functions on $\R$ then their convolution is
$f\ast g(x)=\int_{-\infty}^\infty f(x-t)g(t)\,dt$.  The following proposition gives
the basic properties of convolution.
\begin{prop}\label{c1}
Let $f$ and  $g$  be real-valued functions on $\R$.  
Define $f_x\real$ by $f_x(y) =f(x+y)$ for $x,y\in\R$.
For interval $I=[\alpha,\beta]\subset\R$ and $y\in\R$, 
define $I-y=[\alpha-y,\beta-y]$.

\begin{enumerate}
\item[
{\rm (a)}] 
If $f\ast g$ exists at $x\in\R$ then 
$f\ast g(x) =g\ast f (x)$.
\item[
{\rm (b)}] If $f\in\hk$, $g\in \bv$ and $h\in L^1$ then
$(f\ast g)\ast
h=f\ast(g\ast h)$ on $\R$.
\item[
{\rm (c)}] Let $f\in\hk$.  Suppose that for each compact 
interval $I\subset \R$ there
are constants $K_I$ and $M_I$ such that $|g|\ast |h|(z)\leq K_I$ for all $z\in I$
and the function $y\mapsto h(y)V_{I-y}g$ is in $L^1$.
If $f\ast(g\ast h)$ exists at $x\in\R$ then $(f\ast g)\ast h(x)=f\ast(g\ast h)(x)$.
\item[
{\rm (d)}] $(f\ast g)_x =
f_x\ast g=f\ast g_x$ wherever any one of these convolutions exists.
\item[
{\rm (e)}] ${\rm supp}(f\ast g)\subset \{x+y:x\in{\rm supp}(f), y\in{\rm supp}(g)\}$.
\end{enumerate}
\end{prop}

\noindent
{\bf Proof:}
For (a), (d) and (e), the
$L^1$ proofs hold without change.  See \cite[Proposition 8.6]{folland}.
To prove (b), write 
\begin{eqnarray*}
(f\ast g)\ast h(x) & = & \int_{-\infty}^\infty f\ast g(x-y) h(y)\,dy\\
 & = & \int_{-\infty}^\infty\int_{-\infty}^\infty f(x-y-z)g(z)\,dz\,h(y)\,dy\\
& = & \int_{-\infty}^\infty\int_{-\infty}^\infty f(x-z)g(z-y)h(y)\,dy\,dz\\
& = & f\ast(g\ast h)(x).
\end{eqnarray*}
Lemma~\ref{lemma3}(b) allows us to change the order of $y$ and $z$ integration.
The proof of (c) is similar but now we use
Lemma~\ref{lemma3}(a).  \qed

The next proposition gives some
sufficient conditions for existence of the convolution and some point-wise 
estimates.

\begin{prop}\label{convexist}
\begin{enumerate}
\item[{\rm (a)}] Let $f\in\hk$ and $g\in\bv$.  Then $f\ast g$ exists on $\R$ and
$|f\ast g(x)|\leq \|f\|[\inf|g|+Vg]$ for all $x\in\R$.
\item[{\rm (b)}] Let $f\in\hk_{loc}$ and $g\in\bv$ with the support of $g$
 in the compact interval
$[a,b]$.  Then $f\ast g$ exists on $\R$ and
$|f\ast g(x)|\leq |\int_{x-b}^{x-a}f|\inf_{[a,b]}|g|+\|f\|_{[x-a,x-b]}V_{[a,b]}g$.
\end{enumerate}
\end{prop}

\noindent
{\bf Proof:} 
(a) Using Lemma~\ref{lemma2},
\begin{eqnarray*}
|f\ast g(x)| & = & \left|\int_{-\infty}^\infty f(x-t)g(t)\,dt\right|\\
 & \leq & \left|\int_{-\infty}^\infty f\right|\inf|g|+\|f\|Vg\\
& \leq & \|f\|\left[\inf|g|+Vg\right].
\end{eqnarray*}

(b) Now,
\begin{eqnarray*}
|f\ast g(x)| & = & \left|\int_{a}^b f(x-t)g(t)\,dt\right|\\
 & \leq & \left|\int_{x-b}^{x-a}f\right|\inf_{[a,b]}|g|+\|f\|_{[x-b,x-a]}V_{[a,b]}g.
 \quad\qed
 \end{eqnarray*}

These conditions are sufficient but not necessary for existence of the 
convolution.  Also, 
if $f,g\in\hk$ then $f\ast g$ need not exist at any point.
\begin{example}
{\rm
\begin{enumerate}
\item[(a)]Let $f(x)=\log|x|\sin(x)$ and $g(x)=|x|^{-\alpha}$, 
where $0<\alpha<1$.
Then $f$ and $g$ do not have compact support and 
are not in $\hk$, $\bv$ or $L^p$ ($1\leq p\leq\infty$).
And yet $f\ast g$ exists on $\R$.
\item[(b)]Let $f(x)=\sin(x)/|x|^{1/2}$ and $g(x)=(\sin(x)+\cos(x))/|x|^{1/2}$.
Then $f,g\in\hk$ but $f\ast g$ exists nowhere. \qed
\end{enumerate}
}
\end{example}

When $f\in\hk$ and $g\in L^1\cap\bv$ then $f\ast g$ exists on $\R$ and
we can estimate it in  the Alexiewicz norm.

\begin{prop}\label{c3}
Let $f\in\hk$ and $g\in L^1\cap\bv$. Then $f\ast g$ exists on $\R$
and $\|f\ast g\|\leq \|f\| \|g\|_1$.
\end{prop}

\noindent
{\bf Proof:}
Existence comes from Proposition~\ref{convexist}.  Let $-\infty\leq a<b\leq \infty$.
Using Lemma~\ref{lemma3}(a), we can interchange the repeated integrals, 
\begin{eqnarray}
\int_a^b f\ast g\,dx & = & \int_a^b \int_{-\infty}^\infty f(x-t)g(t)\,dt\,dx\label{conv1}\\
 & = & \int_{-\infty}^\infty g(t)\int_a^b f(x-t)\,dx\,dt.\label{conv2}
\end{eqnarray}
And, 
\begin{eqnarray*}
\left|\int_a^b f\ast g\,dx\right| & \leq & \|g\|_1\lsup_{t\in\R}
\left|\int_a^b f(x-t)\,dx\right|\\
 & = & \|g\|_1\lsup_{t\in\R}\left|\int_{a-t}^{b-t} f\right|\\
 & \leq & \|f\| \|g\|_1.\qed
\end{eqnarray*}

Under suitable conditions on $f$ and $g$, we have the usual interactions between convolution and
Fourier transformation and inversion.

\begin{prop}\label{c4}
If $\fhat$ exists at $s\in\R$ and $g\in L^1\cap\bv$ then 
$\widehat{f\ast g}(s)=\fhat(s)\,\hht{g}(s)$.
\end{prop}

\noindent 
{\bf Proof:}
We have
\begin{eqnarray*}
\widehat{f\ast g}(s) & = & \int_{-\infty}^\infty\eisx  \int_{-\infty}^\infty\left[e^{-ist}f(t)\right]
\left[e^{ist}g(x-t)\right]\,dt\,dx\\
 & = & \int_{-\infty}^\infty e^{-ist}f(t)\int_{-\infty}^\infty e^{-is(x-t)}g(x-t)\,dx\, dt\\
 & = & \fhat(s)\,\hht{g}(s).
\end{eqnarray*}
The interchange of integrals is validated by 
Lemma~\ref{lemma3}(a), since 
\begin{eqnarray*}
\int_{-\infty}^\infty V_{[a,b]}\left[t\mapsto e^{-is(x-t)}g(x-t)\right]\,dx & = &
\int_{-\infty}^\infty V_{[x-b,x-a]}
\left[t\mapsto e^{-ist}g(t)\right]dx\\
 & \leq & 2|s|(b-a)\|g\|_1+2(b-a)Vg. \qed
\end{eqnarray*}

\begin{prop}\label{c5}
If $f$ and $g$ are in $\hk_{loc}$ such that 
$\fhat$ exists almost everywhere,  $\hht{g}\in L^1$, $s\mapsto s\,\hht{g}(s)$ is
in $L^1$ and $\widehat{g}\!\text{{\rm{\LARGE \v{}}}}
=g$ almost everywhere
then
${\displaystyle f\ast g=(\fhat\,\hht{g})\!\!\!\text{{\LARGE \v{}}}}$ on $\R$.
\end{prop}

\noindent
{\bf Proof:}
Let $x\in\R$.  Then $\hht{g}(x-t)$ exists for almost all $t\in\R$.  And,
\begin{eqnarray*}
f\ast g(x) & = & 
\frac{1}{2\pi}\int_{-\infty}^\infty f(t) \int_{-\infty}^\infty e^{is(x-t)}\hht{g}(s)\,ds\,dt\\
 & = & \frac{1}{2\pi}\int_{-\infty}^\infty e^{isx}\,\hht{g}(s) \int_{-\infty}^\infty e^{-ist}f(t)
\,dt\,ds\\
 & = & \fhat\,\hht{g})\!\!\text{{\LARGE \v{}}}\!\!(x).
\end{eqnarray*}
Suppose $\fhat$ exists at $s_0$.  Then 
$
V_I[t\mapsto e^{is(x-t)}e^{is_0t}\,\hht{g}(s)]\leq 2|\hht{g}(s)||s-s_0||I|$ and
the reversal of $s$ and $t$ integration order is by Lemma~\ref{lemma3}(a). \qed

\section{Inversion}
A well-known inversion theorem states that if $\fhat$ and 
$\fhat\text{{\LARGE \v{}}}$ are
in $L^1$ then $f=\fhat\text{{\LARGE \v{}}}$ almost everywhere.  These are rather restrictive
conditions as both $f$ and $\fhat$ must be continuous (almost everywhere)
and vanish at infinity.
In Example~\ref{example1}(a) and (b), $\fhat$ is a multiple of $f$ and $\hht{g}$ is
a multiple of $g$ so we certainly have $f=\fhat\text{{\LARGE \v{}}}$ 
and $g=\widehat{g}\!\text{{\LARGE \v{}}}$ almost everywhere and yet none of these integrals exists
in $L^1$.  However, they do exist in $\hk$.  And, we have a similar inversion
theorem in $\hk$.  First we
need
the following Parseval relation.

\begin{prop}\label{propinv}
Let $\psi$ and $\phi$ be real-valued functions on $\R$.  Suppose $\hht{\psi}$ exists at some
$s_0\in\R$.  Suppose $\phi\in L^1$ and the function $s\mapsto s\phi(s)$ is also in $L^1$.  If
$\int_{-\infty}^\infty \psi\,\hht{\phi}$ exists,  then $\hht{\psi}$ exists almost everywhere and
$\int_{-\infty}^\infty \psi\,\hht{\phi}=\int_{-\infty}^\infty \hht{\psi}\,\phi$.
\end{prop}

\noindent
{\bf Proof:}
Let $f(x)=\psi(x)e^{-is_0x}$ and $g(x,y)=e^{i(s_0-y)x}\phi(y)$.  A simple computation shows
$V_{[a,b]}g(\cdot, y)=O((b-a)y\phi(y))$ as $|y|\to\infty$.  The conditions of 
Lemma~\ref{lemma3}(a)
are satisfied.  
\qed

Now we have the inversion theorem.  The proof uses the method of
summability kernels.  Using
Proposition~\ref{propinv}, one inserts a summability kernel in the
inversion integral.  There is a parameter $z=x+iy$ that is sent to $x_0$, yielding
inversion at $x_0$.  We can actually let $z\to x_0$ in the upper complex plane,
provided the approach is non-tangential.  This is analogous to the 
Fatou theorem
for boundary values of harmonic functions.  Define the upper half
plane by $\pip=\{z=x+iy:x\in\R,\,\,\,y>0\}$.  We identify $\ppip$ with $\R$.  For
$x_0\in\ppip$, we say $z\to x_0$ non-tangentially in $\pip$ if $z\in\pip$ and
$z\to x_0$ such that $|x- x_0|/y\leq C$ for some $C>0$.
\begin{defn}[Summability kernel]\label{summability}
A summability kernel is a function $\Theta\real$ such that
$\Theta\in L^1\cap AC$,
$\Theta(0)=1$,
$s\mapsto s\,\Theta(s)$ is in $L^1$,
$\hht{\Theta}\in L^1\cap\bv$,
$\int_{-\infty}^\infty\hht{\Theta}=2\pi$,
$s\mapsto s\widehat{\Theta}\,'(s)$ is in $L^1$,
and
$x\mapsto V_{[x,\infty)}\hht{\Theta}$ and  $x\mapsto V_{(-\infty, -x]}\hht{\Theta}$ are
$O(1/x)$ as $x\to\infty$.
\end{defn}

\begin{theorem}[Inversion]\label{inversion}
Let $f\real$ such that $\fhat$  exists almost everywhere.  Define $F(x)=\int_{x_0}^xf$ for
$x_0\in\R$.  If $F'(x_0)=f(x_0)$ and 
$f=\fhat\text{{\rm {\LARGE \v{}}}}$
exists at $x_0$ then $f(x_0)=\fhat\text{{\rm{\LARGE \v{}}}}\!\!(x_0)$.
If $\fhat\text{{\rm{\LARGE \v{}}}}$ exists almost everywhere then 
$f=\fhat\text{{\rm{\LARGE \v{}}}}$ almost
everywhere.
\end{theorem}

\noindent
{\bf Proof:}
Let $z=x+iy$ for $x\in\R$ and $y>0$.  Define $\phi_z\real$ by
$\phi_z(s)=\Theta(ys) e^{isx}$, where $\Theta$ is a summability kernel.  
Then $\hht{\phiz}(t)=\hht{\Theta}((t-x)/y)/y$.
And,
\begin{eqnarray}
\frac{1}{2\pi}\int_{-\infty}^\infty \phiz(s)\fhat(s)\,ds  & = & 
\frac{1}{2\pi}\int_{-\infty}^\infty \hht{\phiz}(t)f(t)\,dt\label{inv1}\\
 & = & \frac{1}{2\pi y}\int_{-\infty}^\infty \hht{\Theta}((t-x)/y)f(t)\,dt.
\label{inv2}
\end{eqnarray}
The inversion theorem now follows, provided we can prove the following.
\begin{enumerate}
\item[I.] The conditions of Proposition~\ref{propinv} are satisfied so
that \eqref{inv1} is valid.
\item[II.] As $z\to x_0$ non-tangentially in $\pip$, the left side of
\eqref{inv1} becomes $\fhat\text{{\rm{\LARGE \v{}}}}\!\!(x_0)$.
\item[III.] As $z\to x_0$ non-tangentially in $\pip$, \eqref{inv2} becomes
$f(x_0)$.
\end{enumerate}

I. In Proposition~\ref{propinv}, let
$\psi=f$ and 
$\phi=\phiz$.  We have existence of $\fhat$ at some $s_0\in\R$.
And, $\phiz$ and $s\mapsto s\phiz(s)$ are in $L^1$ if and only if $\Theta$
and $s\mapsto s\,\Theta(s)$ are in $L^1$. Since $\Theta\in L^1$, $\widehat{\phiz}$ is continuous
with limit $0$ at infinity.  So, if $\widehat{\phiz}$ is of bounded variation at infinity,
the integral
$\int_{-\infty}^\infty f\widehat{\phiz}$ will exist.  It suffices to have $\hht{\Theta}$ 
of bounded variation at infinity.  Proposition~\ref{propinv} now applies.

II.  Write the left side of \eqref{inv1} as $(2\pi)^{-1}\int_{-\infty}^\infty 
\left[\Theta(ys)e^{is(x-x_0)}\right]\left[e^{isx_0}\fhat(s)\right]\,ds$.  The function
$s\mapsto e^{isx_0}\fhat(s)$ is in $\hk$.   And, we have $V[s\mapsto
\Theta(ys)e^{is(x-x_0)}] \leq 2V\Theta + 2\|\Theta\|_1 |x-x_0|/y$.  So, for non-tangential
approach, this function is of bounded variation, uniformly as $z\to x_0$.  This allows us
to take the limit inside the integral on the left side of \eqref{inv1}, yielding
$\fhat\text{{\rm{\LARGE \v{}}}}\!\!(x_0)$.

III.  Let $\delta >0$.  Write 
\begin{equation}
\frac{1}{y}\!\lint_{-\infty}^\infty\!\!\! \hht{\Theta}\left(\frac{t-x}{y}\right)f(t)\,dt = 
\frac{1}{y}\!\!\!\lint_{|t-x|<\delta}\!\!\!\hht{\Theta}\left(\frac{t-x}{y}\right)f(t)\,dt  +
\frac{1}{y}\!\!\!\lint_{|t-x|>\delta}\!\!\!\hht{\Theta}\left(\frac{t-x}{y}\right)f(t)\,dt.
\label{inv3}
\end{equation}
Consider the last integral in \eqref{inv3}.  There is $s_0\in\R$ such
that $t\mapsto e^{-is_0t}f(t)$ is in $\hk$.  Now, 
$$
V_{[x+\delta,\infty)}\left[t\mapsto \frac{1}{y}\hht{\Theta}\left(\frac{t-x}{y}\right)e^{is_0t}
\right]\leq \frac{2}{y}V_{[\delta/y,\infty)}\hht{\Theta}+2|s_0|\|\hht{\Theta}\|_1.
$$
With our assumptions on $\hht{\Theta}$, this last expression is bounded as $z\to x_0$.  
And,
when $\Theta\in L^1$, $\Theta\in AC_{loc}$ and $\Theta'\in L^1$ then 
$\hht{\Theta}(t)=o(1/t)$ as $t\to\infty$ \cite[page~20]{benedetto}.  The same
applies on the interval $(-\infty,x-\delta]$.    Hence, taking
the limit $z\to x_0$ inside the integral yields $0$ for each fixed $\delta>0$.

Treat the first integral on the right side of \eqref{inv3} as follows.
Because $\frac{1}{2\pi y}\int_{x-\delta}^{x+\delta}\hht{\Theta}((t-x)/y)\,dt=
\frac{1}{2\pi}\int_{-\delta/y}^{\delta/y}\hht{\Theta}(t)\,dt\to 1$ as $y\to 0^+$, 
we can assume $f(x_0)=0$
(otherwise replace $f(\cdot)$ with $f(\cdot)-f(x_0)$).  
Let $F(t)=\int_{x_0}^t f$.  We have $F'(x_0)=f(x_0)=0$.  And, $\hht{\Theta}(s)=o(1/s)$ as $|s|\to \infty$.
Given $\epsilon >0$, we can
take $0<\delta< 1$ small enough such that $|F(x_0+t)|\leq \epsilon|t|$ and
$|\hht{\Theta}(1/t)|\leq\epsilon|t|$
for all $0< |t|\leq
2\delta$.  Without loss of generality, assume $x\geq x_0$. Take $|z-x_0|\leq \delta$ with
$|x-x_0|/y\leq C$ for some constant $C>0$.  
Integrate by parts,
\begin{eqnarray}
\frac{1}{y}\int_{x-\delta}^{x+\delta}\hht{\Theta}\left(\frac{t-x}{y}\right)f(t)\,dt
 & = & \frac{1}{y}
\left[\hht{\Theta}\left(\frac{\delta}{y}\right)F(x+\delta)-
\hht{\Theta}\left(-\frac{\delta}{y}\right)F(x-\delta)\right]\notag\\
 & & \quad - J_1-J_2-J_3,\label{inv4}
\end{eqnarray}
where $J_1=y^{-2}\int_{x-\delta}^{x_0}\hht{\Theta}'((t-x)/y)F(t)\,dt$,
$J_2=y^{-2}\int_{x_0}^{x}\hht{\Theta}'((t-x)/y)F(t)\,dt$,
$J_3=y^{-2}\int_x^{x+\delta}\hht{\Theta}'((t-x)/y)F(t)\,dt$.
Note that if $0<y\leq\delta^2$ then $y/\delta\leq\delta$ and 
$|x\pm\delta-x_0|< 2\delta$ so
$|\hht{\Theta}(\pm\delta/y)F(x\pm\delta)|/y\leq 2\epsilon^2$.  

Estimate $J_1$ by writing
\begin{eqnarray}
|J_1| & \leq & \frac{1}{y^2}\int_{x-\delta}^{x_0}(x_0-t)\left|\hht{\Theta}'\left(\frac{t-x}{y}\right)
\right|
\left|\frac{F(t)}{x_0-t}\right|\,dt\notag\\
 & \leq & \frac{\epsilon}{y}\int_{-\delta/y}^{(x_0-x)/y}(x_0-x-yt)\left|
\hht{\Theta}'(t)\right|\,dt\notag\\
 & \leq & \epsilon C\, V\hht{\Theta}+\epsilon\int_{-\infty}^0 |t|\left|\hht{\Theta}'(t)\right|\,dt.
\label{inv5}
\end{eqnarray}
Similarly, 
\begin{equation}
|J_3|\leq 
\epsilon C\, V\hht{\Theta}+\epsilon\int_0^\infty t\left|\hht{\Theta}'(
t)\right|\,dt.
\label{inv7}
\end{equation}

For $J_2$ we have
\begin{eqnarray}
|J_2| & \leq & \frac{1}{y}\lsup_{x_0\leq t\leq x}|F(t)|\int_{(x_0-x)/y}^{0}\left|\hht{\Theta}'(t)
\right|
\,dt\notag\\
 & \leq & \epsilon C\, V\hht{\Theta}.\label{inv6}
\end{eqnarray}

Putting \eqref{inv5}, \eqref{inv6} and \eqref{inv7} into \eqref{inv4} now shows that
the first integral on the right side of \eqref{inv3} goes to $0$ as $z\to x_0$ 
non-tangentially.
This completes the proof of part III.  Since, $F'=f$ almost 
everywhere, the proof of the theorem is 
now complete.  
\qed

\begin{remark}
{\rm 
In place of the condition $t\mapsto t\,\hht{\Theta}\,'(t)$ 
is in $L^1$ we can demand that
$\hht{\Theta}$ is increasing on $(-\infty,0)$ and decreasing on $(0,\infty)$.  The
proof of III. then follows with minor changes.  The condition that
$\Theta\in AC$ can also be weakened.
}
\end{remark}
\begin{remark}
{\rm 
The most commonly used summability kernels are 
$$
\begin{array}{lll}
\Theta_1(x)= (1-|x|)\chi_{[-1,1]}(x) &\hht{\Theta}_1(s) = \left[
\frac{\sin(s/2)}{s/2}\right]^2 &\text{ Ces\`{a}ro--Fej\'{e}r}\\
\Theta_2(x)= e^{-|x|} &\hht{\Theta}_2(s) = \frac{2}{1+s^2}
&\text{ Abel--Poisson}\\
\Theta_3(x)= e^{-x^2} &\hht{\Theta}_3(s) = \sqrt{\pi}\, e^{-(s/2)^2}
&\text{ Gauss--Weierstrass}.
\end{array}
$$
The Abel and Gauss kernels are summability kernels according to
Definition~\ref{summability}, while the Ces\`{a}ro kernel does not
satisfy this definition.
}
\end{remark}

\begin{corollary}
Let $f\real$.  Then $f=0$ almost everywhere if and only if
$\fhat=0$ almost everywhere.
\end{corollary}

\noindent
{\bf Proof:}
If $f=0$ almost everywhere then $\fhat=0$ on $\R$.
If $\fhat=0$ almost everywhere then
$\fhat$ exists almost everywhere and
$\fhat\text{{\rm{\LARGE \v{}}}}$
exists  almost everywhere.  Therefore, by the Theorem,  
$\fhat\text{{\rm{\LARGE \v{}}}}=f=0$,
almost everywhere. \qed

Note that the inversion theorem applies to Example~\ref{example1}(a)-(d).  The
condition that $\fhat\text{{\rm{\LARGE \v{}}}}$ exists  almost everywhere cannot be dropped.
The following example shows that existence of $\fhat$ on $\R$ does not guarantee
existence of $\fhat\text{{\rm{\LARGE \v{}}}}$ at any point in $\R$.

\begin{example}\label{example2}
{\rm
Let $f(x)=x^\alpha e^{ix^\nu}$ for $x\geq 0$ and $f(x) =0$ for $x<0$.  Using
the method of Lemma~\ref{lemma1} we see that $\fhat$ exists on $\R$ for
$-1<\alpha<\nu -1$.  And,
\begin{eqnarray}
\fhat(s) & = & \int_{0}^\infty x^{\alpha}e^{i[x^\nu-sx]}\,dx\notag\\
 & = & s^{\tfrac{\alpha+1}{\nu-1}}
 \int_{0}^\infty x^{\alpha}e^{ip[x^\nu-x]}\,dx \quad 
 \left(p=s^{\nu/(\nu-1)}\right).
 \end{eqnarray}
Write
$\phi(x)=x^\nu-x$.  If $\nu>1$ then $\phi$ has a minimum at 
$x_0:=\nu^{-1/(\nu-1)}$.
The method of stationary phase \cite{olver} shows that 
$$
\fhat(s)\sim \sqrt{\tfrac{2\pi}{\nu(\nu-1)}}\,e^{i\pi/4}
x_0^{\alpha-(\nu-2)/2}
e^{i\phi(x_0)s^{\nu/(\nu-1)}}
s^{\frac{2\alpha+2-\nu}{2(\nu-1)}}
$$
as $s\to\infty$.
Let $\nu>2$.
It now follows from Lemma~\ref{lemma1} that when $\nu/2\leq \alpha <\nu-1$,
$\fhat$ exists on $\R$ and $\fhat\text{{\rm{\LARGE \v{}}}}$ diverges at
each point of $\R$.  Note that $f\in\hk$ but neither $f$ nor  $\fhat$
is in any $L^p$ space ($1\leq p\leq \infty$). \qed
}
\end{example}
\section{Appendix}

\begin{lemma}\label{lemma1}
If $\gamma>0$  and $\delta\in\R$ then
\begin{enumerate}
\item[{\rm (a)}]  $\int_0^1 e^{ix^{-\gamma}}x^{\delta}dx$ exists in $\hk$ if
and only if $\gamma+\delta+1>0$.  The integral exists in $L^1$ if and
only if $\delta>-1$.
\item[{\rm (b)}] $\int_1^\infty e^{ix^{\gamma}}x^{\delta}dx$ exists in $\hk$ if
and only if $\gamma>\delta+1$.  The integral exists in $L^1$ if and
only if $\delta<-1$.
\end{enumerate}
\end{lemma}

\noindent
{\bf Proof:} In (a), integrate by parts to get
$$
\int_0^1 e^{ix^{-\gamma}}x^{\delta}dx=
\frac{i}{\gamma}\left[e^i-\llim_{x\to 0^+}e^{ix^{-\gamma}}x^{\gamma+\delta+1}\right]-\frac{i(
\gamma+\delta+1)}{\gamma}\int_0^1 e^{ix^{-\gamma}}x^{\gamma+\delta}dx.
$$
The limit exists if and only if $\gamma+\delta+1>0$, the last integral then converging
absolutely.  Case (b) is similar.  For $L^1$ convergence, we simply take the absolute value
of each integrand. \qed

\begin{lemma}\label{lemma2}
Let $[a,b]\subset\rbar$ and let $f\in\hk_{[a,b]}$ and $g\in\bv_{[a,b]}$.
Then
$$
\left|\int_a^bfg\right|\leq 
\left|\int_a^bf\right|
\linf_{[a,b]}|g|+
\|f\|_{[a,b]}V_{[a,b]}g.
$$
\end{lemma}

\noindent
{\bf Proof:} Given
$\epsilon>0$, take $c\in [a,b]$ such that $|g(c)|\leq \epsilon + \inf_{[a,b]}|g|$.
Integrate by parts:
\begin{eqnarray*}
\int_a^bfg & = & \int_a^cfg + \int_c^bfg\\
 & = & g(c)\int_a^b f-\int_a^c\left(\int_a^xf\right)dg(x) + \int_c^b\left(\int_x^bf\right)dg(x).
 \end{eqnarray*}
 And, 
 \begin{eqnarray*}
 \left|\int_a^bfg\right| & \!\leq\! & \left[\epsilon + \linf_{[a,b]}|g|
 \right]\left|\int_a^b f\right| + \lsup_{a\leq x\leq c}\left|\int_a^xf\right|V_{[a,c]}g
 + \lsup_{c\leq x\leq b}\left|\int_x^bf\right|V_{[c,b]}g\\
  & \!\leq\! & \left[\epsilon + \linf_{[a,b]}|g|
   \right]\left|\int_a^b f\right|
 + \|f\|_{[a,b]}V_{[a,b]}g. \qed
   \end{eqnarray*}

This lemma is an extension of inequalities proved in \cite{rieszlivingston}
and \cite{celidze} (Theorem~ 45, page~36).  Changing $g$ on a set of measure $0$, such as
a singleton, does not affect the integral of $fg$ but  can make the infimum of $|g|$ equal to
zero.  However, this reduction in $\inf|g|$ is reflected by a corresponding increase in 
$Vg$.  This redundancy can be eliminated by replacing $g$ with its normalised version, i.e.,
for each $x\in[a,b)$ replace $g(x)$ with $\lim_{t\to x+}g(t)$ and redefine $g(b)=0$.
Then the inequality becomes $\left|\int_a^bfg\right|\leq
\|f\|_{[a,b]}V_{[a,b]}g$.

The following lemma on interchange of iterated integrals
is an extension of Theorem~57 on page~58 of \cite{celidze}.
\begin{lemma}\label{lemma3}
Let $f\in\hk$ and let $g\fn\R^2\to\R$.  Let $\meas$ denote the measurable subsets of $\R$.
For each $(A,B)\in\bv\times\meas$, define the iterated integrals
\begin{eqnarray*}
I_1(A,B) & = & \int_{x\in A}\int_{y\in B}f(x)g(x,y)\,dy\,dx\\
I_2(A,B) & = & \int_{y\in B}\int_{x\in A}f(x)g(x,y)\,dx\,dy.\\
\end{eqnarray*}
\begin{enumerate}
\item[{\rm (a)}]
Assume that for each compact interval $I\subset \R$ there are constants
$M_I>0$ and $K_I>0$ such that
$\int_{\R}V_Ig(\cdot, y)\,dy\leq M_I$ and, for all $x\in I$,
$\|g(x,\cdot)\|_1\leq K_I$. 
If $I_1$ exists on $\R\times\R$ then $I_2$ exists on
$\bv\times\meas$ and $I_1=I_2$ on $\bv\times\meas$.
\item[{\rm (b)}]
Assume there exist $M,G\in L^1$ such that, for almost all $y\in\R$, $Vg(\cdot, y)\leq M(y)$ and, 
for all $x\in\R$, $|g(x,y)|\leq G(y)$.
Then $I_1=I_2$ on $\bv\times\meas$.
\end{enumerate}
\end{lemma}

\noindent
{\bf Proof:} (a) Let ${\cal I}$ be the open intervals in $\R$.  First prove $I_1=I_2$ on 
${\cal I}\times{\cal I}$.  Fix $(a,b)$ and $(\alpha,\beta)$ in $I$.  
For $-\infty<a<t<\infty$, define
\begin{equation}
H_a(t)=I_2((a,t),(\alpha,\beta))=\int_{\alpha}^\beta\int_a^tf(x)g(x,y)\,dx\,dy.\label{app2}
\end{equation}
We will establish the equality of $I_1$ and $I_2$ by appealing to the necessary and sufficient
conditions for interchanging repeated integrals \cite[Corollary 6]{talviladiff}.  For this,
we need to show that $H_a$ is in $ACG_*$ and that we can differentiate under the integral sign
in \eqref{app2}.
Let $F(x)=\int_{-\infty}^x f$.  Integrate by parts,
\begin{equation}
H_a(t)=\left[F(t)-F(a)\right]
\int_\alpha^\beta g(t,y)\,dy-\int_{\alpha}^\beta\int_a^t\left[F(x)-F(a)\right]d_1g(x,y)dy.\label{app1}
\end{equation}
The integrator of the Riemann--Stieltjes integral over $x\in[a,t]$ is denoted $d_1g(x,y)$.  Now, by Lemma~\ref{lemma2},
$|H_a(t)|  \leq  \|f\|[K_{[a,b]}+M_{[a,b]}]$, and $I_2(A,B)$ exists for all $A,B\in{\cal I}$
with $A$ bounded.

We have $F\in ACG_*(\R)$.  So, there
are $E_n\subset\R$ such that 
$\R=\cup E_n$ and $F$ is $AC_*$ on each $E_n$, i.e., for each
$n\geq 1$, given $\epsilon>0$, there is $\delta >0$ such that if $(s_i,t_i)$ are disjoint with
$s_i,t_i\in E_n$ and $\sum|s_i-t_i|<\delta$ then $\sum\|f\|_{(s_i,t_i)}<\epsilon$.  Fix
$n\geq 1$
with $E_n$, $\epsilon$ and $\delta$ as above.  Suppose $(\sigma_i, \tau_i)$ are disjoint with
$\sigma_i,\tau_i\in E_n$ and $\sum|\sigma_i-\tau_i|<\delta$.  With no loss of generality, we
may assume $E_n$ is a subset of a compact interval $[c,d]$.  Then 
$$
\lsup_{[p,q]\subset[\sigma_i,\tau_i]}\left|H_a(p)-H_a(q)\right|\leq
\|f\|_{[\sigma_i,\tau_i]}\left[K_{[c,d]}+M_{[c,d]}\right].
$$
It follows that $H_a\in ACG_*(\R)$.

Now show that we can differentiate under the integral sign to compute $H_a'(t)$.
Let $0<|h|<1$ and $t\in\R$ such that $F'(t)=f(t)$.  Then
\begin{eqnarray*}
\left|\frac{1}{h}\int_t^{t+h}f(x)g(x,y)\,dx\right| & \leq & 
\lsup_{0<|h|<1}\left|\frac{F(t+h)-F(t)}{h}\right||g(t,y)|\\
 & & \quad +\lsup_{0<|h|<1}\left|\frac{1}{h}\right|\|f\|_{[t-|h|,t+|h|]}V_{[t-1,t+1]}g(\cdot,y).
\end{eqnarray*}
It now follows from dominated convergence that $H_a'(t)=f(t)\int_\alpha^\beta g(t,y)\,dy$ for
almost all $t\in\R$.  And, by \cite[Corollary 6]{talviladiff}, $I_1(A,B)=I_2(A,B)$ for all
$A,B\in {\cal I}$ with $A$ bounded.

By assumption, $I_1(\R,\R)$ exists.  For $a\in\R$,
\begin{eqnarray*}
\int_a^\infty f(x)\int_\alpha^\beta g(x,y)\,dy\,dx & = & \llim_{t\to\infty}
\int_a^t f(x)\int_\alpha^\beta g(x,y)\,dy\,dx\\
 & = &  \llim_{t\to\infty} \int_\alpha^\beta\int_a^t f(x)g(x,y)\,dx\,dy\\
 & = & \llim_{t\to\infty} H_a(t).
\end{eqnarray*}
Similarly, $\llim_{t\to-\infty} H_a(t)$ exists.  Therefore, $H_{-\infty}$ is continuous on ${\overline \R}$
and hence in $ACG_*({\overline \R})$.  It follows from \cite[Corollary 6]{talviladiff} that
$I_1(A,B)=I_2(A,B)$ for all
$A,B\in {\cal I}$.

We have equality of $I_1$ and $I_2$ on $\bv\times\meas$ upon replacing $f$ with $f\chi_A$ and $g(x,\cdot)$
with $g(x,\cdot)\chi_B$ where $A\in \bv$ and $B\in \meas$.

(b) This is similar to part (a), but now the conditions on  $g$ ensure 
the existence
of $I_2$ on $\R\times\R$.  As in (a), $H_a\in ACG_*(\R)$.  To show $H_a$ is continuous on ${\overline \R}$,
note that
$ |\int_\alpha^\beta\int_a^tf(x)g(x,y)\,dx\,dy|\leq \|f\|(\|G\|_1+\|M\|_1)$ and $\lim_{t\to\infty}
\int_a^tf(x)g(x,y)\,dx$ exists for almost all $y\in\R$.  Whence,
$\lim_{t\to\infty}H_a(t)$ exists and $H_{-\infty}$ is continuous on ${\overline \R}$.  Using
\cite[Corollary 6]{talviladiff}, we now have equality of $I_1$ and $I_2$ on $\R\times\R$ and
hence on $\bv\times\meas$.\qed

\end{document}